\documentclass[12pt]{article}
\usepackage{amssymb,amsthm,amsmath,amsfonts,latexsym,tikz,hyperref,ytableau,authblk}
\usetikzlibrary {arrows.meta}
\usepackage[hmargin=.9in,vmargin=1in]{geometry}

\newtheorem{thm}{Theorem}[section]
\newtheorem{prop}[thm]{Proposition}
\newtheorem{cor}[thm]{Corollary}
\newtheorem{lem}[thm]{Lemma}
\newtheorem{conj}[thm]{Conjecture}
\newtheorem{exa}[thm]{Example}

\DeclareMathOperator{\fusc}{fusc}
\DeclareMathOperator{\CW}{CW}
\DeclareMathOperator{\depth}{dp}
\DeclareMathOperator{\height}{ht}

\DeclareMathOperator{\irr}{Irr}
\DeclareMathOperator{\rgf}{rgf}
\DeclareMathOperator{\SL}{SL}
\newcommand{\hb}{\ol{h}}
\newcommand{\fs}{\mathsf{s}}
\newcommand{\fst}{\tilde{\mathsf{s}}}

\newcommand{\doubleleftbrace}{\{\!\{}
\newcommand{\doublerightbrace}{\}\!\}}

\newcommand{\ben}{\begin{enumerate}}
\newcommand{\een}{\end{enumerate}}
\newcommand{\ble}{\begin{lem}}
\newcommand{\ele}{\end{lem}}
\newcommand{\bth}{\begin{thm}}
\renewcommand{\eth}{\end{thm}}
\newcommand{\bpr}{\begin{prop}}
\newcommand{\epr}{\end{prop}}
\newcommand{\bco}{\begin{cor}}
\newcommand{\eco}{\end{cor}}
\newcommand{\bcon}{\begin{conj}}
\newcommand{\econ}{\end{conj}}
\newcommand{\bde}{\begin{defn}}
\newcommand{\ede}{\end{defn}}
\newcommand{\bex}{\begin{exa}}
\newcommand{\eex}{\end{exa}}
\newcommand{\barr}{\begin{array}}
\newcommand{\earr}{\end{array}}
\newcommand{\btab}{\begin{tabular}}
\newcommand{\etab}{\end{tabular}}
\newcommand{\beq}{\begin{equation}}
\newcommand{\eeq}{\end{equation}}
\newcommand{\bea}{\begin{eqnarray*}}
\newcommand{\eea}{\end{eqnarray*}}
\newcommand{\bal}{\begin{align*}}
\newcommand{\bce}{\begin{center}}
\newcommand{\ece}{\end{center}}
\newcommand{\bpi}{\begin{picture}}
\newcommand{\epi}{\end{picture}}
\newcommand{\bpp}{\begin{picture}}
\newcommand{\epp}{\end{picture}}
\newcommand{\bfi}{\begin{figure} \begin{center}}
\newcommand{\efi}{\end{center} \end{figure}}
\newcommand{\bprf}{\begin{proof}}
\newcommand{\eprf}{\end{proof}\medskip}

\newcommand{\bsl}{\begin{slide}{}}
\newcommand{\esl}{\end{slide}}
\newcommand{\bfr}{\begin{frame}}
\newcommand{\efr}{\end{frame}}

\newcommand{\eqed}[1]{$\textcolor{white}{\qed}\hfill{\dil#1}\hfill\qed$}

\newcommand{\ol}{\overline}

\newcommand{\hs}[1]{\hspace{#1}}
\newcommand{\hso}[1]{\hspace{-1pt}}
\newcommand{\vs}[1]{\vspace{#1}}
\newcommand{\qmq}[1]{\quad\mbox{#1}\quad}

\newcommand{\emp}{\emptyset}

\newcommand{\sbe}{\subseteq}

\newcommand{\iso}{\cong}

\newcommand{\zh}{\hat{0}}
\newcommand{\oh}{\hat{1}}

\newcommand{\ptn}{\vdash}

\newcommand{\lt}{\lhd}
\newcommand{\gt}{\rhd}

\newcommand{\case}[4]{\left\{\barr{ll}#1&\mbox{#2}\\#3&\mbox{#4}\earr\right.}

\def\<{\langle}
\def\>{\rangle}

\newcommand{\ra}{\rightarrow}

\newcommand{\be}{\beta}

\newcommand{\ep}{\epsilon}

\newcommand{\la}{\lambda}

\newcommand{\La}{\Lambda}

\newcommand{\1}{{\bf 1}}

\newcommand{\bbN}{{\mathbb N}}

\newcommand{\bbQ}{{\mathbb Q}}
\newcommand{\bbR}{{\mathbb R}}

\newcommand{\cD}{{\cal D}}

\newcommand{\cF}{{\cal F}}
\newcommand{\cG}{{\cal G}}
\newcommand{\cH}{{\cal H}}

\newcommand{\cJ}{{\cal J}}

\newcommand{\dil}{\displaystyle}

\begin{document}
\pagestyle{plain}

\title{Hyperbinary partitions and 
$q$-deformed rationals}
\author[1]{Thomas McConville}
\author[2]{James Propp}
\author[3]{Bruce E. Sagan}
\affil[1]{Department of Mathematics, Kennesaw State University, Kennesaw, GA 30144}
\affil[2]{Department of Mathematics and Statistics, University of Massachusetts Lowell, Lowell, MA 01854}
\affil[3]{Department of Mathematics, Michigan State University, East Lansing, MI 48824}

\date{\today\\[10pt]
	\begin{flushleft}
	                                       \\[5pt]
	\end{flushleft}}

\maketitle

\begin{abstract}
A hyperbinary partition of the nonnegative integer $n$ is a partition where every part is a power of $2$ and every part appears at most twice. 
We give three applications of the length generating function for such partitions, denoted by $h_q(n)$.
Morier-Genoud and Ovsienko defined the $q$-analogue of a rational number $[r/s]_q$ in various ways,
most of which depend directly or indirectly
on the continued fraction expansion of $r/s$.
As our first application we show that $[r/s]_q=q\,h_q(n-1)/h_q(n)$ where $r/s$ occurs as the $n$th entry in the Calkin-Wilf enumeration of the non-negative rationals.
Next we consider fence posets which are those which can be obtained from a sequence of chains by alternately pasting together maxima and minima.  For every $n$ we show there is a fence poset $\cF(n)$ whose lattice of order ideals is isomorphic to the poset of hyperbinary partitions of $n$ ordered by refinement.
For our last application, Morier-Genoud and Ovsienko also showed that $[r/s]_q$ can be computed by taking products of certain matrices which are $q$-analogues of the standard generators for the special linear group $\SL(2,\bbR)$.  We express the entries of these products in terms of the polynomials $h_q(n)$.
\end{abstract}

\section{Introduction}

There has been much work in recent years on 
Stern’s diatomic sequence (e.g.~\cite{CW:rr}), 
fence posets (e.g.~\cite{OR:rpfpu}), and 
$q$-deformed rational numbers (e.g.~\cite{MGO:qdr1}), 
with links between these topics. 
We strengthen these links by bringing into the foreground {\em hyperbinary partitions}. 
These are partitions in which all parts are powers of two 
and in which no part appears more than twice. 
These have appeared in the literature on Stern’s diatomic sequence, 
but it has not been noticed that these objects relate to order ideals in fence posets 
and that a natural statistic on these partitions gives 
a nice way to construct the $q$-deformed rational numbers, 
avoiding explicit reliance on continued fractions.
We explain those additional links.

In view of the centrol role to be played by hyperbinary partitions, we first establish some definitions and notation about integer partitions in general.
If $\la$ is an integer partition then we will write it either as a weakly decreasing sequence of integers $\la=(\la_1,\la_2,\ldots,\la_\ell)$ or in terms of multiplicities
$$
\la=\doubleleftbrace 1^{m_1}, 2^{m_2}, \ldots, n^{m_n}\doublerightbrace 
$$
where
$$
m_i = m_i(\la) =\text{the number of $i$'s in $\la$}.
$$
When using multiplicity notation in examples, we will often dispense with the commas and multiset braces. 
When the multiplicity $m_i$ is 1, we write $i^1$ as $i$; 
when the multiplicity $m_i$ is 0, we omit $i^0$ entirely.
For example, the integer partition $(4,1,1)$ can be written as $1^2 4$. 
We may also choose to list the parts in an order other than increasing, writing $1^2 4$ as $4 1^2$ or even $141$.
Regardless of the notation chosen, if $\la$ is a partition of $n$ (meaning that the sum of its parts is $n$) then we will write $\la\ptn n$.
The {\em length} of $\la$ is
$$
\ell(\la) = \text{the number of parts of $\la$} = 
\sum_{i} m_i (\la).
$$
Returning to our example, $\ell(4,1,1)=3$.

Call a partition $\eta$ {\em hyperbinary} if
\ben
\item each part  is a power of $2$, and
\item the multiplicity of each part is at most $2$.
\een
It appears that Wilf coined this term.
The first in-depth study of such partitions was made by 
Reznick~\cite{rez:sbpf}, though antecedents can be found going as far back as Stern~\cite{ste:uzf}.
Let
\beq
\label{H(n)}
H(n) = \{\eta \mid \text{$\eta$ is a hyperbinary partition of $n$}\}
\eeq
and
$$
h(n) = \#H(n)
$$
where we will use $\#S$ or $|S|$ for the cardinality of a set $S$.
For example,
$$
H(10) = \{82,\ 81^2,\ 4^2 2,\ 4^2 1^2,\ 4 2^2 1^2\}
$$
so that
$$
h(10) = 5.
$$
We introduce the generating function
$$
h_q(n) = \sum_{\eta\in H(n)} q^{\ell(\eta)}.
$$
For instance, 
$$
h_q(10) = q^2 + 2q^3 + q^4 + q^5.
$$
Clearly $h_1(n) = h(n)$.   We will give three applications using $h_q(n)$.

Our first application, which is in the next section, involves the Calkin-Wilf sequence $\CW(n)$, $n\ge0$.  This sequence is defined as the ratio $\CW(n)=\fusc(n)/\fusc(n+1)$ where $\fusc(n)$ is Stern's diatomic sequence as reinvented by Dijkstra (see \eqref{fuscDef}).  The Calkin-Wilf sequence goes through each nonnegative rational number exactly once.  Mourier-Genoud and Ovsienko gave a way of associating with any rational number $r/s$ a $q$-analogue which is a rational function $[r/s]_q$.  Our main result of this section is that one can calculation the $q$-analogue of $\CW(n)$ using the polynomials $h_q(n)$.  More precisely, we show in Theorem~\ref{thm:qrat} that
$$
[\CW(n)]_q = q\, \frac{h_q(n-1)}{h_q(n)}.
$$

In Section~\ref{php}, we consider the poset (partially ordered set) $\cH(n)$ of hyperbinary partitions of $n$ under the refinement ordering.  A fence is a poset obtained by taking a sequence of chains and alternately identifying their maxima and minima.  Our principal result here is the isomorphism in Theorem~\ref{thm:mainbij} which shows that
$\cH(n)\iso \cJ(\cF(n))$ where $\cF(n)$ is the fence associated with $n$, and $\cJ(P)$ is the distributive lattice of all lower order ideals of the poset $P$ under inclusion.

Section~\ref{m} is devoted to the study of certain $q$-analogues of the standard generators of $\SL(2,\bbR)$, see~\eqref{LR}.  Morier-Genoud and Ovsienko showed that their rational $q$-analogues can be computed using certain products of these matrices.  We prove in Theorem~\ref{thm:MnEnt} that the entries of such products can be easily computed using the $h_q(n)$.

We end with a section devoted to open questions and avenues for future research.

\section{A $q$-analogue of the Calkin-Wilf sequence}

Let $\bbN$ and $\bbQ$ be the nonnegative integers and the rationals, respectively.
{\em Stern's diatomic sequence}, also known as the {\em Stern-Brocot sequence} or the {\em obfuscating sequence}, can be defined inductively as 
$\fusc(0) = 0$, $\fusc(1) = 1$, and for $n\ge 1$,
\beq
\label{fuscDef}
\begin{array}{rcl}
\fusc(2n) &=& \fusc(n),\\
\fusc(2n+1)&=& \fusc(n+1)+\fusc(n)
\end{array}
\eeq
(see Table 1).
To our knowledge, Stern~\cite{ste:uzf} was the first person to study this sequence.   
The $\fusc$ notation was coined by Dijkstra~\cite[pp.\ 215-216]{dij:swc}.  
For a history of this sequence, see the article of Northshield~\cite{nor:sds}.

The {\em Calkin-Wilf sequence} is defined for all $n\ge0$ by
$$
\CW(n) = \frac{\fusc(n)}{\fusc(n+1)}.
$$
This function has the property that 
for each rational number $r/s \geq 0$
there is a unique integer $n \geq 0$ satisfying $\CW(n) = r/s$.
Calkin and Wilf introduced this sequence in~\cite{CW:rr} and related the $\fusc$ function to hyperbinary partitions.

\begin{table}
\centering
\begin{tabular}{rclcl}
$n$ & $\fusc_n$ & $\CW_n$ & $\fusc_n(q)$ & $\CW_n(q)$\\[1ex]
\hline
 0 & 0 & 0 & 0 & 0 \\[1ex]
 1 & 1 & 1 & 1 & $\frac{1}{q}$ \\[1ex]
 2 & 1 & $\frac12$ & $q$ & $\frac{1}{1+q}$ \\[1ex]
 3 & 2 & 2 & $q + q^2$ & $\frac{1+q}{q}$ \\[1ex]
 4 & 1 & $\frac13$ & $q^2$ & $\frac{q}{1+q+q^2}$   \\[1ex]
 5 & 3 & $\frac32$ & $q + q^2 + q^3$ & $\frac{1 + q + q^2}{q + q^2}$ \\[1ex]
 6 & 2 & $\frac23$ & $q^2 + q^3$ & $\frac{1 + q}{1 + q + q^2}$ \\[1ex]
 7 & 3 & 3 & $q^2 + q^3 + q^4$ & $\frac{1 + q + q^2}{q}$ \\[1ex]
 8 & 1 & $\frac14$ & $q^3$ & $\frac{q^2}{1 + q + q^2 + q^3}$ \\[1ex]
 9 & 4 & $\frac43$ & $q + q^2 + q^3 + q^4$ & $\frac{1 + q + q^2 + q^3}{q + q^2 + q^3}$ \\[1ex]
10 & 3 & $\frac35$ & $q^2 + q^3 + q^4$ & $\frac{1 + q + q^2}{1 + 2 q + q^2 + q^3}$ \\[1ex]
11 & 5 & $\frac52$ & $q^2+2q^3+q^4+q^5$ & $\frac{1+2q+q^2+q^3}{q+q^2}$ 
\end{tabular}
\caption{The functions $\fusc_n$, $\CW_n$, $\fusc_n(q)$ and $\CW_n(q)$}
\end{table}

We mention here a method for computing $n$ from $r/s$ 
that essentially is described in~\cite{CW:bp} 
and deserves to be better known.
Recall that every positive rational number $r/s$
has two representations as continued fractions,
that is, representations of the form
\beq
\label{ConFra}
\frac{r}{s}= 
a_1+\frac{1}{\dil{a_2 + \frac{1}{\ddots + \dil{\frac{1}{a_m}}}}}
\eeq
where $a_1 \geq 0$ and $a_2,\dots,a_m \geq 1$;
for instance, $7/3$ can be written as both
$2+1/3$ (with $m=2$) and as $2+1/(2+1/1)$ (with $m=3$).
Given $r/s$, pick the representation with odd length. Create a binary string consisting of
$a_1$ 1's followed by $a_2$ 0's followed by $a_3$ 1's
followed by \dots followed by $a_m$ 1's.
Reverse it and one obtains the binary representation 
of the unique $n$ satisfying $\CW(n) = r/s$. 
For instance, with $r/s = 7/3 = 2+1/(2+1/1)$ 
we form the bit-string 11001 whose reversal $10011$ is the binary expansion of the number 19,
and one can check that $\fusc(19) = 7$ and $\fusc(20) = 3$
yielding $\CW(19) = 7/3$.

We will need three operations on partitions.  Suppose $\la=\doubleleftbrace 1^{m_1}, 2^{m_2}, \ldots, n^{m_n}\doublerightbrace $ and $\mu=\doubleleftbrace 1^{k_1}, 2^{k_2}, \ldots, n^{k_n}\doublerightbrace $.  Then their {\em sum} is the partition
\beq
\label{la+mu}
\la+\mu = \{\!\{ 1^{m_1+k_1}, 2^{m_2+k_2}, \ldots, n^{m_n+k_n}\}\!\}.
\eeq
If $k_i\le m_i$ for all $i$ then their {\em difference} is
\beq
\label{la-mu}
\la-\mu = \doubleleftbrace 1^{m_1-k_1}, 2^{m_2-k_2}, \ldots, n^{m_n-k_n}\doublerightbrace .
\eeq
If $t$ is a positive rational number and 
$\la=(\la_1,\la_2,\ldots,\la_\ell)$ then their {\em product} is
\beq
\label{tla}
t\la = (t\la_1,t\la_2,\ldots,t\la_\ell)
\eeq
provided that $t\la_i\in\bbN$ for all $i$. 
We extend these operations to sets $\La$ of partitions by letting
\begin{align}
\La+\mu & =\{\la+\mu \mid \la\in\La\},  \label{La+mu}\\
\La-\mu & =\{\la-\mu \mid \la\in\La\}, \label{La-mu}\\
t\La & =\{t\la \mid \la\in\La\}, \label{tLa}
\end{align}
provided the sets of the right sides of the equal signs are sets of partitions.

With respect to the three operations, we have
\begin{align}
\ell(\la+\mu) &= \ell(\la)+\ell(\mu),\label{ell+}\\
\ell(\la-\mu) &= \ell(\la)-\ell(\mu),\label{ell-}\\
\ell(t\la) &=\ell(\la).\label{ellr}
\end{align}

We now show that the sets $H(n)$ defined by~\eqref{H(n)} have a nice recursive structure.  Let $\ep$ denote the empty partition and $\uplus$ denote the disjoint-union operation on sets.  The following result is in~\cite{CW:rr}, but we include its proof for completeness.
\bpr[\cite{CW:rr}]
\label{Hrr}
We have $H(-1)=\emp$, $H(0) = \{ \ep\}$, and for $n\ge1$
\begin{align}
H(2n-1) &= 2H(n-1)+(1),\label{H(2n-1)}\\
H(2n) &= 2H(n) \uplus [2H(n-1) +(1^2)].\label{H(2n)}
\end{align}
\epr
\bprf
For equation~\eqref{H(2n-1)}, note that if $\eta\in H(2n-1)$ then $m_1(\eta)=1$ since $\eta$ is a hyperbinary partition of an odd number.  Thus $\eta-(1)$ is a hyperbinary partition of $2n-2$ with all parts at least $2$.  It follows that $\eta-(1)=2\psi$ for some $\psi\in H(n-1)$ and the desired equality follows.

Now consider~\eqref{H(2n)}.  If $\eta\in H(2n)$ then $1$ appears with multiplicity zero or two.  In the first case  $\eta=2\psi$ where $\psi\in H(n)$.  In the second,
$\eta-(1^2) = 2\chi$ where $\chi\in H(n-1)$.  This finishes the proof of the equality and of the proposition.
\eprf

We now show that $h_q(n-1)$ can be used as a $q$-analogue of $\fusc(n)$.
\bpr
\label{hqrr}
We have $h_q(-1)=0$, $h_q(0) = 1$, and for $n\ge1$
\begin{align}
h_q(2n-1) &= q h_q(n-1),\label{hq(2n-1)}\\
h_q(2n) &= h_q(n) + q^2 h_q(n-1).\label{hq(2n)}
\end{align}
\epr
\bprf
In view of the properties of the length function (equations~\eqref{ell+}, \eqref{ell-}, and~\eqref{ellr}), this result is just a translation of  Proposition~\ref{Hrr} into the language of generating functions.
\eprf

Comparison of the previous proposition with the definition of  the Stern sequence in~\eqref{fuscDef} prompts the following definition.  Define the 
{\em $q$-Stern sequence} to be the polynomial sequence where $\fusc_q(0) = 0$ and for $n\ge 1$,
$$
\fusc_q(n) = h_q(n-1).
$$
Translating the previous proposition into the language of the $\fusc_q$ polynomials gives $\fusc_q(0) = 0$, $\fusc_q(1) = 1$, and
\beq
\label{fuscqRr}
\begin{array}{rcl}
\fusc_q(2n) &=& q \fusc_q(n),\\
\fusc_q(2n+1)&=& \fusc_q(n+1) + q^2 \fusc_q(n)
\end{array}
\eeq
for $n \geq 1$.
Similarly, we define the {\em $q$-Calkin-Wilf sequence} to be the sequence of rational functions
$$
\CW_q(n) = \frac{\fusc_q(n)}{\fusc_q(n+1)} = \frac{h_q(n-1)}{h_q(n)}
$$
for $n \geq 1$, with $\CW_q(0) = 0$.

There is another way to obtain a related $q$-analogue of the Calkin-Wilf sequence.  Morier-Genoud and Ovsienko~\cite{MGO:qdr1,MGO:qdr2,MGO:qdr3} 
found a way to associate with every rational number $r/s\in\bbQ$ a rational function $[r/s]_q\in\bbQ(q)$ which has many interesting properties and connections to various branches of mathematics.   Suppose that $r/s$ is a positive rational number and 
consider the continued fraction expansion of $r/s$ as in~\eqref{ConFra}.  The notation for this expansion is $r/s=[a_1,a_2,\ldots,a_m]$.
Now define the {\em $q$-analogue} of $r/s$, $[r/s]_q$, to be the rational function obtained by taking the continued fraction for $r$ and making the replacements
$$
a_i \text{ becomes } 
\begin{cases}
[a_i]_q & \text{if $i$ is odd,}\\
[a_i]_{q^{-1}} & \text{if $i$ is even,}
\end{cases}
$$
and
$$
\text{the $1$ in the $i$th numerator becomes } 
\begin{cases}
q^{a_i} & \text{if $i$ is odd,}\\
q^{-a_i} & \text{if $i$ is even,}
\end{cases}
$$
where $[a_i]_q$ denotes the ordinary $q$-integer
$1 + q + q^2 + \cdots + q^{a_i - 1}$.
The result of these substitutions is denoted $[r/s]_q = [a_1,a_2,\ldots,a_m]_q$ and the initial part of the fraction is
$$
\left[\frac{r}{s}\right]_q= 
[a_1]_q+\frac{q^{a_1}}{\dil{[a_2]_{q^{-1}} + \frac{q^{-a_2}}{\ddots}}}.
$$
It is easy to see that $[r/s]_q$ does not depend on which of the two continued fraction expansions one starts with.

Now one could ask if there is a relationship between  $\CW_q(n)$ and the $q$-analogue given by  
$$
[\CW(n)]_q=\left[\frac{\fusc(n)}{\fusc(n+1)}\right]_q.
$$
To see what the relationship is, we will need the fact, proved by Morier-Genoud and Ovsienko,
 that for all rational numbers $r$ we have
\beq
\label{r/s+1}
\left[\frac{r}{s} + 1\right]_q = q\left[\frac{r}{s}\right]_q + 1.
\eeq
\bth
\label{thm:qrat}
For all $n\ge0$ we have
$$
[\CW(n)]_q =  q \CW_q(n).
$$
\eth
\bprf
We induct on $n$ where, as we will usually do, the base case will be omitted because it is easy.  We first consider odd arguments $n$.  Then, using the recurrence relations~\eqref{fuscqRr}, we obtain
$$
\CW_q(2n+1) = \frac{\fusc_q(2n+1)}{\fusc_q(2n+2)} = \frac{q^2 \fusc_q(n)+ \fusc_q(n+1)}{q \fusc_q(n+1)} = q\ \CW_q(n) + \frac{1}{q}.
$$
Thus, by induction and~\eqref{r/s+1},
$$
q\CW_q(2n+1) = q^2 \CW_q(n) + 1 = q(q\CW_q(n))+1
= q([\CW(n)]_q)+1 = [\CW(n)+1]_q,
$$
On the other hand, 
$$[\CW(2n+1)]_q =  \left[\frac{\fusc(2n+1)}{\fusc(2n+2)} \right]_q=  \left[\frac{\fusc(n)+ \fusc(n+1)}{\fusc(n+1)}\right]_q = [\CW(n) + 1]_q 
$$
Comparing the expressions for $q\CW_q(2n+1)$ and $[\CW(2n+1)]_q $ completes this case.

As far as even arguments go, 
\begin{align}
q\CW_q(2n) & = \frac{q\fusc_q(2n)}{\fusc_q(2n+1)}\notag \\[5pt]
&= \frac{q^2 \fusc_q(n)}{q^2 \fusc_q(n)+ \fusc_q(n+1)} \notag \\[5pt]
&=\frac{q}{q + \frac{\fusc_q(n+1)}{q\fusc_q(n)}} \notag \\[5pt]
&=\frac{q}{q + \frac{1}{q\CW_q(n)}}.\label{CW_q(2n)}
\end{align}
Similarly,
$$
[\CW(2n)]_q =
\left[\frac{\fusc(2n)}{\fusc(2n+1)} \right]_q =  \left[\frac{\fusc(n)}{\fusc(n)+\fusc(n+1)}\right]_q =
\left[\frac{1}{1+\frac{\fusc(n+1)}{\fusc(n)}}\right]_q
$$
yielding
\beq
\label{CW(2n)_q}
[\CW(2n)]_q = \left[\frac{1}{1+\frac{1}{\CW(n)}}\right]_q.
\eeq
Now there are two subcases depending on whether $\CW(n)\ge 1$ or $\CW(n)<1$.  We will do the former as the latter is similar.

Suppose $\CW(n)=[a_1,a_2,\ldots,a_m]$. Then, since $\CW(n)\ge 1$ we have that 
$$\frac{1}{1 + \frac{1}{\CW(n)}} = [0,1,a_1,a_2,\ldots,a_m].$$
Combining this with~\eqref{CW(2n)_q} 
and the definition of a rational $q$-analogue gives
$$
[\CW(2n)]_q 
= [0]_q+\frac{q^0}{\dil{[1]_{q^{-1}} + \frac{q^{-1}}{[\CW(n)]_q}}}
=\frac{1}{\dil{1 + \frac{q^{-1}}{[\CW(n)]_q}}}
=\frac{q}{\dil{q + \frac{1}{[\CW(n)]_q}}}.
$$
Comparing this expression with~\eqref{CW_q(2n)} and using the induction hypothesis completes the proof of the theorem.
\eprf

\section{The poset of hyperbinary partitions of $n$}
\label{php}

\begin{figure}
    \centering
\begin{tikzpicture}
\filldraw(0,0) circle(.1);
\filldraw(0,1) circle(.1);
\filldraw(-1,2) circle(.1);
\filldraw(1,2) circle(.1);
\filldraw(0,3) circle(.1);
\draw (0,0)--(0,1)--(-1,2)--(0,3)  (0,1)--(1,2)--(0,3);
\draw(0,-.5) node{$\doubleleftbrace 4,2^2,1^2\doublerightbrace $};
\draw(1,1) node{$\doubleleftbrace 4^2,1^2\doublerightbrace $};
\draw(-2,2) node{$\doubleleftbrace 8,1^2\doublerightbrace $};
\draw(2,2) node{$\doubleleftbrace 4^2,2\doublerightbrace $};
\draw(0,3.5) node{$\doubleleftbrace 8,2\doublerightbrace $};
\draw(0,-1.5) node{$\cH(10)$};
\draw(4,2) node{$\iso$};
\begin{scope}[shift={(7.5,0)}]
\filldraw(0,0) circle(.1);
\filldraw(0,1) circle(.1);
\filldraw(-1,2) circle(.1);
\filldraw(1,2) circle(.1);
\filldraw(0,3) circle(.1);
\draw (0,0)--(0,1)--(-1,2)--(0,3)  (0,1)--(1,2)--(0,3);
\draw(0,-.5) node{$0122$};
\draw(1,1) node{$0202$};
\draw(-2,2) node{$1002$};
\draw(2,2) node{$0210$};
\draw(0,3.5) node{$1010$};
\draw(0,-1.5) node{$\cD(10)$};
\end{scope}
\end{tikzpicture}    
    \caption{The posets $\cH(10)$ and $\cD(10)$ \label{cHcD10}}
\end{figure}

Let $\cH(n)$ denote the poset of hyperbinary partitions under the {\em refinement partial order},
where we say $\mu$ {\em refines} $\la$
(in symbols, $\mu \leq \la$)
if the parts of $\la$ can be subdivided to produce the parts of $\mu$.
An equivalent way to state this definition is that 
 the parts of $\mu$ can be grouped together
so that, adding the parts in each group, one obtains the parts of $\la$.
For example, $\cH(10)$ is displayed on the left in Figure~\ref{cHcD10}.
The poset $\cH(n)$ was studied by Brunetti and D'Aniello~\cite{BDA:gch} who used it to study how the length of a hyperbinary expansion of $n$ (see the definition of such an expansion in the next paragraph) is related to $n$ itself.
Our aim is to show that $\cH(n)$ is isomorphic to the lattice of ideals of a corresponding fence poset.  For any undefined terms used from the theory of partially ordered sets, see the texts of Sagan~\cite{sag:aoc} or Stanley~\cite{sta:ec1}.
It's worth mentioning that the poset of {\em all} partitions of $n$ is {\em not} a lattice under refinement order when $n\ge5$; for instance, 
the partitions $41$ and $32$ both cover the partitions $311$ and $221$
so the former two do not have a meet (coarsest common refinement)
while the latter two do not have a join (finest common coarsening).

It will be convenient to use hyperbinary expansions rather than hyperbinary partitions.  
Suppose that the binary expansion of $n$ is
$$
\be(n) := b_1 b_2\ldots b_k,
$$
in other words
$$
n = b_1 2^{k-1} + b_2 2^{k-2} + \ldots + b_k.
$$
Note our nonstandard convention of having $b_1$ be the coefficient of the highest power of $2$, $b_2$ for the next-highest, and so forth.  This will make the indexing simpler when we describe the isomorphism.  
A {\em hyperbinary expansion} of $n$ is 
$$
d = d_1 d_2 \ldots d_k
$$
having the same length as the binary expansion $\be(n)$ where $d_i\in\{0,1,2\}$ for all $i$ and
$$
n =  d_1 2^{k-1} + d_2 2^{k-2} + \ldots + d_k.
$$
Note that there may be some initial zeros in a hyperbinary expansion forced by the fact that it has the same number of digits as the binary expansion.  
For example, if $n=10$ then the largest power of $2$ in its binary expansion is $2^3$ so all hyperbinary expansions must have length $3+1=4$. 
More specifically, $d=0122$ is a hyperbinary expansion for $10$ since it has length $4$ and
$$
10= 0\cdot 2^3 + 1\cdot 2^2 + 2\cdot 2^1 + 2.
$$
Given a sequence $d=d_1\ldots d_k$ of zeros, ones, and twos, we let 
\begin{align*}
s(d) &= \text{ the integer for which $d$ is a hyperbinary expansion}\\
&=\sum_{i=1}^k d_i 2^{k-i}.
\end{align*}
Note that we may need to adjust the number of initial zeros to make the length of $d$ correct.  So, as just noted, 
$s(0122)=10$.
For a more refined invariant, 
we let 
$$
s_i(d) = s(d_1\ldots d_i).
$$
For example, $s_3(10210)=s(102)=1\cdot 2^2+0\cdot 2+2\cdot 1=6$. 

There is a clear bijection between hyperbinary partitions $\eta$ of $n$ and hyperbinary expansions $d$ of $n$ obtained by mapping $\eta$  to $d=d_1\ldots d_k$, where $2^{k-1}$ is the largest power of $2$ in $\be(n)$ and $d_i$ is the multiplicity of $2^{k-i}$ in $\eta$.
Thus the set $\cD(n)$ of hyperbinary expansions of $n$ inherits a poset structure induced by $\cH(n)$.  See Figure~\ref{cHcD10} for this isomorphism when $n=10$.

The following lemma will be useful.  
It shows that our definition of $\cH(n)$ coincides with that in~\cite{BDA:gch}.
We write 
$x\lt y$ if $x$ is {\em covered by} $y$, i.e., $x<y$ and there is no $z$ with $x<z<y$.
\ble
\label{cov}
Element $d=d_1\ldots d_k\in\cD(n)$ covers exactly the elements which can be obtained from $d$ by replacing some adjacent pair $d_i 0$ where $d_i>0$ with the pair
$(d_i-1)2$.
\ele
\bprf
In $\cH(n)$ the partial order is refinement.  So a partition $\eta$ covers those partitions which can be formed from it by replacing a part $2^j$ with two parts $2^{j-1}+2^{j-1}$.  Note that by the hyperbinary restriction, this can only be done if there are no parts of the form $2^{j-1}$ already in $\eta$.  Translating in terms of hyperbinary expansions, these are the covers described in the lemma.

To show that these are the only ones, suppose that $d=d_1\ldots d_k$ covers 
$c=c_1\ldots c_k$.  Then $c$ is obtained by refining a single part of $d$, since if two or more parts were refined then refining only one of them would give an element strictly between the two.  The possible refinements of a part $2^j$ as a hyperbinary partition are all of the form 
$$
2^j = 2^{j-1} + 2^{j-2} + \cdots + 2^{l+1} + 2^l + 2^l
$$
for some $l<j$.  Let $d_r d_{r+1}\ldots d_s$ be the corresponding digits in $d$ with $d_r\ge1$ parts equal to $2^j$
(so $r=k-j$ and $s=k-l$).  Thus in $c$ we have 
$$
c_r c_{r+1} \ldots c_s = (d_r-1) (d_{r+1}+1) (d_{r+2}+1)\ldots (d_{s-1}+1)(d_s+2).
$$
In order for this to be a valid  hyperbinary expression, we must have $d_s=0$ and $d_i=0$ or $1$ for all $r<i<s$.  For $r<t<s$, let  $d_t$ be the rightmost $1$. (If all these $d_i$ are zero then a similar argument works using $t=r$.)  Replace $d_t d_{t+1} \ldots d_s=10 \ldots 0$ with 
$01\ldots1 2$.  The resulting $d'$ satisfies $d'<d$.  Now iterate this process, starting with the rightmost $1$ in the factor $d_r\ldots d_{t-1}0$ of $d'$.  This will produce a sequence $d>d'>\ldots>d'' = c$ which shows that $d$ did not cover $c$ to begin with.  This contradiction ends the proof.
\eprf

A poset $P$ has a {\em maximum} if there is an element $\oh$ such that $\oh\ge x$ for all $x\in P$.  Dually, a {\em minimum} is $\zh$ satisfying $\zh\le x$ for all $x\in P$.  The next proposition can also be found in~\cite{BDA:gch}, but we include a proof for completeness.
\bpr\label{prop:bounds}
We have the following.
\ben
\item[(a)] Poset $\cD(n)$ has a maximum, denoted $\oh(n)$, which is the binary expansion of $n$.
\item[(b)] Poset $\cD(n)$ has a minimum, denoted $\zh(n)$, which is the unique hyperbinary expansion whose zeros form a prefix of $\zh(n)$.  
\een
\epr
\bprf
(a) 
Let $d=d_1\ldots d_k\in\cD(n)$. 
Suppose $d$ has at least one entry equal to $2$, and choose $i$ to be the minimum index where $d_i=2$.
Since $n<2^k$ we have $i>1$.
By Lemma~\ref{cov}, $d$ is covered by the element obtained by replacing $d_{i-1}2$ with $(d_{i-1}+1)0$. 

Hence, any maximal element of $\cD(n)$ only has $0$'s and $1$'s. 
The only such element is the binary expansion of $n$.

\medskip

(b) 
Since $\cD(n)$ is finite, it has minimal elements (those which do not cover any other element).  And from Lemma~\ref{cov} it is clear that any minimal element has the form specified in the proposition.  So it suffices to prove that there exists a unique minimal element.

Suppose, to the contrary the $c=c_1\ldots c_k$ and $d=d_1\ldots d_k$ are both minimal in $\cD(n)$.  Let $i$ be the leftmost index in which they differ.  Without loss of generality, suppose $c_i<d_i$.   We will show that $s(c)<s(d)$ so that they cannot both be in $\cD(n)$.  Since $c_1\ldots c_{i-1}=d_1\ldots d_{i-1}$, we need only consider the contribution of $c_i\ldots c_k$ and $d_i\ldots d_k$ to $s(c)$ and $s(d)$, respectively.  But, since $c_i<d_i\le 2$ the largest possible value of $s(c)$ is when $c':=c_i\ldots c_k = 1 2\ldots 2$.  Also, by the placement of zeros in $d$, the smallest value of $s(d)$ with $c_i<d_i$ is when $d':=d_i\ldots d_k=2 1\ldots 1$.  But, from the definition of the function $s$, we have $s(c')=2^{k-i+1}+2^{k-i}-2$ while 
$s(d')=2^{k-i+1}+2^{k-i}-1$.  So $s(c)<s(d)$ as desired.
\eprf

For the next result, we need another concept.
Again consider the binary expansion $\be(n)=b_1 b_2\ldots b_k$.   The {\em principal prefix} of $\be(n)$ is
\beq
\label{pp}
p(\be(n)) = b_1 b_2\ldots b_r
\eeq
where $b_{r+1}$ is the rightmost $0$ in $\be(n)$.  
For the rest of this section we will use $r$ for the length of the principal prefix.
Note  that if $b_i=1$ for all $i$ then, because there is no such zero, $p(\be(n))=\emp$ (the empty string).  
For example, if $n=75$ then $\be(75) = 1001011$ and $p(\be(75))=1001$. 

\bco\label{cor:min}
If $n=2^k-1$, then $\zh(n)=\oh(n)=1^k$.
Else, if $p(\be(n))=b_1\ldots b_r$ then
$$
\zh(n)=0(b_2+1)\ldots(b_r+1)21^{k-r-1}.
$$

\eco

\bprf
If $n=2^k-1$, then $1^k$ is the unique hyperbinary expansion of $n$.

Suppose $n\ne 2^k-1$, and let $c=0(b_2+1)\ldots(b_r+1)21^{k-r-1}$.
Since the binary expansion of $n$ only has $0$'s and $1$'s, the entry $b_i+1$ is either $1$ or $2$.
Hence, $c$ only has one zero entry at the beginning, so by Proposition~\ref{prop:bounds}~(b), it remains to show that this word is a hyperbinary expansion of $n$.

Recall $b_1=1$, $b_{r+1}=0$, and $b_i=1$ for $r+2\le i\le k$. Thus,
\begin{align*}
    s(c) &= \left(\sum_{i=2}^r (b_i+1)2^{k-i}\right) + 2\cdot 2^{k-r-1} + \left(\sum_{i=r+2}^k 1\cdot 2^{k-i}\right)\\
    &= 2^{k-r} + \sum_{i=2}^r 2^{k-i} + \sum_{i=2}^k b_i\cdot 2^{k-i}\\
    &= 2^{k-1} + \sum_{i=2}^k b_i\cdot 2^{k-i}\\
    &= n
\end{align*}
as desired.
\eprf

The following lemma will be used to compare two partial orderings on $\cD(n)$.

\ble\label{altord}
Suppose $\le$ and $\preceq$ are partial orders on the same finite set $P$. Assume that for all $x,y\in P$, if $x\preceq y$, then either
\begin{itemize}
\item $x=y$,
\item there exists $z$ such that $x<z\preceq y$, or
\item there exists $w$ such that $x\preceq w<y$. 
\end{itemize}
Then for all $x,y\in P$, if $x\preceq y$, then $x\le y$.
\ele

\bprf
With respect to the partial order $\le$, we define the \emph{depth} of an element $x$ to be the length of the longest chain of $(P,\le)$ whose minimum element is $x$. The \emph{height} of $x$ is the length of the longest chain of $(P,\le)$ whose maximum element is $x$. Throughout this proof, we only consider depth and height with respect to $\le$ rather than $\preceq$. Let $\depth(x)$ and $\height(x)$ denote the depth and height of $x$, respectively.

To prove the lemma, we proceed by induction on $\depth(x)+\height(y)$.
For the base case, consider elements $x,y\in P$ such that $\depth(x)=0=\height(y)$.  Now suppose $x\preceq y$ so that one (or more) of the three conditions in the statement of the lemma must hold. If $x=y$ then $x\le y$, and we are done. It now suffices to prove that the other two conditions  are impossible. 
If there is $z\in P$ such that $x<z\preceq y$, then $\depth(x)>\depth(z)\ge 0$
which is a contradiction to the base case assumption. Likewise, if $w\in P$ such that $x\preceq w<y$ then $\height(y)>\height(w)\ge 0$. 

Now let $k\ge 1$, and suppose the lemma holds for any $x,y\in P$ such that $\depth(x)+\height(y)<k$. Let $x,y\in P$ such that $\depth(x)+\height(y)=k$ and $x\preceq y$.
Again, one of the three conditions of the lemma must hold and the proof breaks up into cases depending on them.

If $x=y$, then $x\le y$, as desired.

For the second case, suppose there exists $z$ such that $x<z\preceq y$. Then $\depth(z)<\depth(x)$, which implies $\depth(z)+\height(y)<k$. Hence, $z\le y$ by the inductive hypothesis. By transitivity, we deduce $x\le y$.

For the third case, suppose there exists $w$ such that $x\preceq w<y$. Then $\height(w)<\height(y)$, which implies $\depth(x)+\height(w)<k$. Similarly to the second case, we have $x\le w$ by the inductive hypothesis.  So, again, $x\le y$.
\eprf

In the next proposition, we give an alternate interpretation of the partial order on hyperbinary expansions of $n$. 

\bpr\label{prop:partsumscompare}
Suppose $c=c_1\ldots c_k$ and $d=d_1\ldots d_k$ are in $\cD(n)$.  Then $c\le d$ if and only if for all $1\le i\le k$ we have $s_i(c)\le s_i(d)$.
\epr

\bprf
For the forward direction, it suffices to show that if $c\lt d$ then the inequalities hold.  From Lemma~\ref{cov}, we have that $c$ is obtained from $d$ by replacing a pair $d_j0$ where $d_j>0$ with $(d_j-1)2$.  It follows that $s_j(d) = s_j(c)+1$, and $s_i(d)=s_i(c)$ for all 
$i\neq j$.  

For the reverse implication, suppose $c,d\in \cD(n)$ such that $s_i(c)\le s_i(d)$ for all $i$. If $s_i(c)=s_i(d)$ for all $i$, then $c=d$, and we are done. Otherwise, let $p$ be the smallest index such that $s_p(c)<s_p(d)$. By the minimality of $p$, we have $c_i=d_i$ for $i<p$ and $c_p<d_p$. Hence, $s_p(d)-s_p(c)=d_p-c_p$.

By Lemma~\ref{altord}, it is enough to show that there is a cover $c\lt e$ such that $s_i(e)\le s_i(d)$ for all $i$, or there is a cover $f\lt d$ such that $s_i(c)\le s_i(f)$ for all $i$.

From the equality
\[ 
s_p(c)\cdot 2^{k-p} + \sum_{i=p+1}^k c_i\cdot 2^{k-i} = n = s_p(d)\cdot 2^{k-p} + \sum_{i=p+1}^k d_i\cdot 2^{k-i},
\]
we have
\begin{align*}
    (d_p-c_p)\cdot 2^{k-p} &= (s_p(d) - s_p(c))\cdot 2^{k-p}\\ 
    &= \sum_{i=p+1}^k (c_i-d_i)\cdot 2^{k-i}\\
    &\le (c_{p+1}-d_{p+1})\cdot 2^{k-p-1} + \sum_{i=p+2}^k2\cdot 2^{k-i}\\
    &= (c_{p+1}-d_{p+1})\cdot 2^{k-p-1} + 2\cdot(2^{k-p-1}-1)
\end{align*}
From this and the fact that $c_p<d_p$ we obtain
\[ 
(c_{p+1}-d_{p+1})\cdot 2^{k-p-1} \ge 2^{k-p} - (2^{k-p}-2)>0,
\]
which implies $c_{p+1}>d_{p+1}$. Hence, either $c_{p+1}=2$ or $d_{p+1}=0$, or both.

Consider the case $c_{p+1}=2$. Since $c_p<d_p\le 2$, there is a cover $c\lt e$ where $e$ is obtained from $c$ by replacing $c_p2$ with $(c_p+1)0$. In this case, we have 
$s_p(e)=s_p(c)+1\le s_p(d)$.  And if $i\neq p$, then $s_i(e)=s_i(c)\le s_i(d)$.
So $s_i(e)\le s_i(d)$ for all $i$ as desired.

Finally, consider the case $d_{p+1}=0$. Also, $d_p>c_p\ge0$.  So there is a cover $f\lt d$ where $f$ is obtained from $d$ by replacing $d_p0$ with $(d_p-1)2$. In this case, $s_p(f)=s_p(d)-1\ge s_p(c)$. And if $i\neq p$, then $s_i(f)=s_i(d)\ge s_i(c)$.
So, again, the desired conclusion holds.
\eprf

We will now get fence posets involved.  
The {\em $n$th fence}, $\cF(n)$, is the poset constructed from the principal prefix $p(\be(n))=b_1 b_2 \ldots b_r$ as follows.  The elements of $\cF(n)$ will be $x_1,x_2,\ldots,x_r$.  Covers will only be between adjacent elements in this list, where
we start with the element $x_1$ and inductively define for $i\ge2$
$$
\case{x_i\lt x_{i-1}}{if $b_i=0$,}{x_i\gt x_{i-1}}{if $b_i=1$.}
$$
As an example, suppose $n=75$.  Recalling that  $p(\be(75))=1001$, we obtain
$$
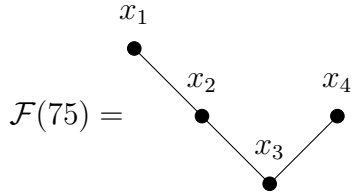
\begin{tikzpicture}[scale=.9]
\draw(-1,1) node{$\cF(75)=$};
\filldraw(0,2) circle(.1);
\filldraw(1,1) circle(.1);
\filldraw(2,0) circle(.1);
\filldraw(3,1) circle(.1);
\draw (0,2)--(2,0)--(3,1);
\draw(0,2.5) node{$x_1$};
\draw(1,1.5) node{$x_2$};
\draw(2,0.5) node{$x_3$};
\draw(3,1.5) node{$x_4$};
\end{tikzpicture}
$$
where the two ``down" covers from $x_1$ to $x_2$ and from $x_2$ to $x_3$ come from the two zeros of $1001$ while the ``up" cover from $x_3$ to $x_4$ comes from the final $1$.

\begin{figure}
    \centering
\begin{tikzpicture}[scale=.9]
\filldraw(0,0) circle(.1);
\filldraw(0,1) circle(.1);
\filldraw(-1,2) circle(.1);
\filldraw(1,2) circle(.1);
\filldraw(0,3) circle(.1);
\draw (0,0)--(0,1)--(-1,2)--(0,3)  (0,1)--(1,2)--(0,3);
\draw(0,-.5) node{$\emp$};
\draw(1,1) node{$\{x_2\}$};
\draw(-2,2) node{$\{x_1,x_2\}$};
\draw(2,2) node{$\{x_2,x_3\}$};
\draw(0,3.5) node{$\{x_1,x_2,x_3\}$};
\draw(0,-1.5) node{$\cJ(\cF(10))$};
\end{tikzpicture}    
    \caption{The poset $\cJ(\cF(10))$}
    \label{cJcA10}
\end{figure}

Let $\cJ(P)$ be the distributive lattice of all lower order ideals of $P$ ordered by containment. 
As an example, consider $\cJ(\cF(10))$.  Now $\be(10)=1010$ so that $p(\be(10))=101$ and
$$
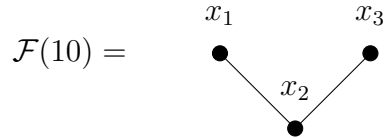
\begin{tikzpicture}
\draw(-1,1) node{$\cF(10)=$};
\filldraw(1,1) circle(.1);
\filldraw(2,0) circle(.1);
\filldraw(3,1) circle(.1);
\draw (1,1)--(2,0)--(3,1);
\draw(1,1.5) node{$x_1$};
\draw(2,0.5) node{$x_2$};
\draw(3,1.5) node{$x_3$};
\end{tikzpicture}
$$
is the correponding poset.  The lattice of order ideals $\cJ(\cF(10))$ is displayed in Figure~\ref{cJcA10}.

The set $\bbN^r$ is partially ordered such that for $u,v\in\bbN^r$, we have $u\leq v$ if and only if $u_i\leq v_i$ for all $i$.
This poset is a distributive lattice where the meet and join operations may be explicitly defined as
\begin{align*}
    u\wedge v &= (\min(u_1,v_1),\ldots,\min(u_k,v_k)),\ \mbox{and}\\
    u\vee v &= (\max(u_1,v_1),\ldots,\max(u_k,v_k)).    
\end{align*}

We construct an isomorphism between $\cD(n)$ and $\cJ(\cF(n))$ by identifying each poset with a sublattice of $\bbN^r$.

If $I$ is a subset of $\cF(n)$, its indicator function is 
\[ \chi_I( i ) = \begin{cases} 0\ &\mbox{if }x_i\notin I\\1 \ &\mbox{if }x_i\in I\end{cases}. \]
We will identify the indicator function $\chi_I$ with its sequence of values $(\chi_I(1),\ldots,\chi_I(r)) \in\bbN^r$.
It is straight-forward to check that 
the function
$\chi(I)= \chi_I$ is a lattice embedding of $\cJ(\cF(n))$ into $\{0,1\}^r$.
When $n=10$, the embedding is illustrated in Figure~\ref{fig:iso}.

\begin{figure}
    \centering
\begin{tikzpicture}
\filldraw(0,0) circle(.1);
\filldraw(0,1) circle(.1);
\filldraw(-1,2) circle(.1);
\filldraw(1,2) circle(.1);
\filldraw(0,3) circle(.1);
\draw (0,0)--(0,1)--(-1,2)--(0,3)  (0,1)--(1,2)--(0,3);
\draw(0,-.5) node{$\emp$};
\draw(1,1) node{$\{x_2\}$};
\draw(-2,2) node{$\{x_1,x_2\}$};
\draw(2,2) node{$\{x_2,x_3\}$};
\draw(0,3.5) node{$\{x_1,x_2,x_3\}$};
\draw(0,-1.5) node{$\cJ(\cF(10))$};
\draw[|->] (3,1.5)--(4.6,1.5);
\draw(3.8,2) node{$\chi$};
\begin{scope}[shift={(7.5,0)}]
\filldraw(0,0) circle(.1);
\filldraw(0,1) circle(.1);
\filldraw(-1,2) circle(.1);
\filldraw(1,2) circle(.1);
\filldraw(0,3) circle(.1);
\draw (0,0)--(0,1)--(-1,2)--(0,3)  (0,1)--(1,2)--(0,3);
\draw(0,-.5) node{$(0,0,0)$};
\draw(1,1) node{$(0,1,0)$};
\draw(-2,2) node{$(1,1,0)$};
\draw(2,2) node{$(0,1,1)$};
\draw(0,3.5) node{$(1,1,1)$};
\draw(0,-1.5) node{$\chi(\cJ(\cF(10)))$};
\end{scope}
\end{tikzpicture}    
    \caption{Embedding  $\cJ(\cF(10))$ in $\{0,1\}^3$
    \label{fig:iso}}
\end{figure}

It remains to show that $\cD(n)$ is isomorphic to the same sublattice of $\{0,1\}^r$ as $\cJ(\cF(n))$.
Given $c\in\cD(n)$, let $\mathsf{s}(c) = (s_1(c),\ldots,s_r(c))$, where 
\[ 
s_j(c)=s(c_1\cdots c_j) = \sum_{i=1}^j c_i\cdot 2^{j-i}
\]
for $j\in[k]$.
\bpr
\label{sEmb}
The map $c\mapsto\mathsf{s}(c)$ embeds $\cD(n)$ as a subposet of $\bbN^r$.
\epr

\bprf
For $c,d\in \cD(n)$, we have $c\le d$ if and only if $s_i(c)\le s_i(d)$ for all $1\le i\le k$ by Proposition~\ref{prop:partsumscompare}.
Hence, we have a poset embedding of $\cD(n)$ into $\bbN^k$.
By definition of the principal prefix and Lemma~\ref{cov}, we have $c_i=d_i=1$ for $i>r+1$.
So, $s_i(c)=s_i(d)$ whenever $i>r$ and the map $c\mapsto \mathsf{s}(c)$ is a poset embedding of $\cD(n)$ into $\bbN^r$.
\eprf

For example, consider $n=10$. 
The binary expansion $\beta(10)=1010$ has principal part $101$, so $r=3$. 
We compute 
\[ 
\mathsf{s}(0210) = (s(0),\ s(02),\ s(021)) = (0,2,5).
\]
Applying $\mathsf{s}$ to each hyperbinary expansion of $10$ gives the poset embedding $\cD(10)\ra \bbN^3$ in the first line of Figure~\ref{smap}.

By a direct calculation, we have the following useful identity.

\ble\label{lem:recursum}
For $c\in\cD(n)$ and $1< i \le k$, we have $s_i(c) = 2\cdot s_{i-1}(c) + c_i$.
\ele

\bpr\label{prop:coversum}
If $c\lt d$ is a cover in $\cD(n)$, then $\mathsf{s}(c)\lt \mathsf{s}(d)$ is a cover in $\bbN^r$.
\epr

\bprf
Suppose $c\lt d$ is a cover in $\cD(n)$.
By Lemma~\ref{cov}, there is an index $j$ such that $c$ is obtained from $d$ by replacing $d_j0$ with $(d_j-1)2$.
It is clear that $s_i(c) = s_i(d)$ for $i<j$.
We compute
\begin{align*}
s_j(c) &= 2\cdot s_{j-1}(c) + c_j\\
&= 2\cdot s_{j-1}(d) + (d_j - 1) = s_j(d) - 1,    
\end{align*}
and
\begin{align*}
s_{j+1}(c) &= 4\cdot s_{j-1}(c) + 2\cdot c_j + 2\\
&= 4\cdot s_{j-1}(d) + 2\cdot d_j = s_{j+1}(d).    
\end{align*}
Since $c_i=d_i$ for $i>j+1$, we deduce through Lemma~\ref{lem:recursum} that $s_i(c)=s_i(j)$ for $i>j+1$.
Therefore, $\mathsf{s}(c)$ is covered by $\mathsf{s}(d)$ in $\bbN^r$.
\eprf

An {\em order filter} of a poset $P$ is a subposet $F$ such that $x\in F$ and $y \ge x$ implies $y\in F$. An order filter $F$ is {\em principal} if it is generated by a single element, i.e. there exists $x\in P$ such that $F = \{y\in P: y\ge x\}$.

For any $u\in\bbN^r$, the poset $\bbN^r$ is isomorphic to the principal order filter $F$ generated by $u$ via the map $F\ra \bbN^r$ where $v\mapsto v-u$.
Setting $\zh=\zh_{\cD(n)}$, we define for any $c\in\cD(n)$ the sequence 
$$
\tilde{\mathsf{s}}(c) = \mathsf{s}(c) - \mathsf{s}(\zh).
$$
Continuing our example when $n=10$, the second line of Figure~\ref{smap} shows the effect of $\fst$ on $\cD(10)$ and illustrates the next proposition.

\begin{figure}
    \centering
\begin{tikzpicture}[scale=.8]
\filldraw(0,0) circle(.1);
\filldraw(0,1) circle(.1);
\filldraw(-1,2) circle(.1);
\filldraw(1,2) circle(.1);
\filldraw(0,3) circle(.1);
\draw (0,0)--(0,1)--(-1,2)--(0,3)  (0,1)--(1,2)--(0,3);
\draw(0,-.5) node{$0122$};
\draw(1,1) node{$0202$};
\draw(-2,2) node{$1002$};
\draw(2,2) node{$0210$};
\draw(0,3.5) node{$1010$};
\draw(0,-1.5) node{$\cD(10)$};
\draw[|->] (3,1.5)--(4.6,1.5);
\draw(3.8,2) node{$\fs$};
%
\begin{scope}[shift={(7.5,0)}]
\filldraw(0,0) circle(.1);
\filldraw(0,1) circle(.1);
\filldraw(-1,2) circle(.1);
\filldraw(1,2) circle(.1);
\filldraw(0,3) circle(.1);
\draw (0,0)--(0,1)--(-1,2)--(0,3)  (0,1)--(1,2)--(0,3);
\draw(0,-.5) node{$(0,1,4)$};
\draw(1,1) node{$(0,2,4)$};
\draw(-2,2) node{$(1,2,4)$};
\draw(2,2) node{$(0,2,5)$};
\draw(0,3.5) node{$(1,2,5)$};
\draw(0,-1.5) node{$\mathsf{s}(\cD(10))$};
\end{scope}
\end{tikzpicture}

\vs{20pt}
\begin{tikzpicture}[scale=.8]
\draw(0,0) node{$\rule{2in}{0pt}$};
\draw[|->] (3,1.5)--(4.6,1.5);
\draw(3.8,2) node{$\fst$};
\begin{scope}[shift={(7.5,0)}]
\filldraw(0,0) circle(.1);
\filldraw(0,1) circle(.1);
\filldraw(-1,2) circle(.1);
\filldraw(1,2) circle(.1);
\filldraw(0,3) circle(.1);
\draw (0,0)--(0,1)--(-1,2)--(0,3)  (0,1)--(1,2)--(0,3);
\draw(0,-.5) node{$(0,0,0)$};
\draw(1,1) node{$(0,1,0)$};
\draw(-2,2) node{$(1,1,0)$};
\draw(2,2) node{$(0,1,1)$};
\draw(0,3.5) node{$(1,1,1)$};
\draw(0,-1.5) node{$\fst(\cD(10))$};
\end{scope}
\end{tikzpicture} 
    \caption{Embedding $\cD(10)$ in $\bbN^3$ and $\{0,1\}^3$ 
    \label{smap}}
\end{figure}




\bpr
The map $c\mapsto\tilde{\mathsf{s}}(c)$ embeds $\cD(n)$ as a subposet of $\{0,1\}^r$ such that $\tilde{\mathsf{s}}(\zh)=(0,\ldots,0)$ and $\tilde{\mathsf{s}}(\oh)=(1,\ldots,1)$.
\epr

\bprf
By definition, we have $\tilde{\mathsf{s}}(\zh)=(0,\ldots,0)$. 
By the discussion above and Proposition~\ref{sEmb}, the map $\tilde{\mathsf{s}}$ gives an embedding of $\cD(n)$ into $\bbN^r$ that sends the minimum element of $\cD(n)$ to the minimum element of $\bbN^r$.   If $n=2^k-1$ then $r=0$ and the proposition is trivial, so we assume $n\neq 2^k-1$.

By Corollary~\ref{cor:min},   the $j$-th entry of $\tilde{\mathsf{s}}(\oh)$ where $j\in[r]$ is 
\begin{align*}
    \fst(\oh)_j&=s_j(\oh) - s_j(\zh) \\
    &= \left(\sum_{i=1}^j b_i\cdot 2^{j-i}\right) - \left(\sum_{i=2}^j (b_i+1)\cdot 2^{j-i}\right)\\
    &= 2^{j-1} - 2^{j-2} - \cdots - 1\\
    &= 1.
\end{align*}
Since  $\fst(\zh)\le\fst(c)\le\fst(\oh)$ for all $c\in\cD(n)$, we have that the embedding is into 
$\{0,1\}^r$.
\eprf

Given a lattice $L$, a subposet $S$ is a {\em join subsemilattice} if for all $x,y\in S$, the join $x\vee_L y$ is in $S$.
In this case, the subposet $S$ is a join semilattice where $x\vee_S y=x\vee_L y$ for any $x,y\in S$.
A {\em meet subsemilattice} is defined dually.
A {\em sublattice} of $L$ is both a join subsemilattice and meet subsemilattice.

The sublattice property can be detected by the following local lemma, which can be deduced from \cite[Lemma 9-2.10]{rea:ltpr}. 

\ble\label{lem:sublattice_local}
Let $L$ be a finite lattice, and let $S$ be a bounded subposet of $L$. 
If for all $a,b,c\in S$ such that $a\lt_S b$ and $a\lt_S c$ the join $b\vee_L c$ is in $S$, then $S$ is a join subsemilattice of $L$.
Dually, if for all $b,c,d\in S$ such that $b\lt_S d$ and $c\lt_S d$ the meet $b\wedge_L c$ is in $S$, then $S$ is a meet subsemilattice of $L$.
\ele

\bpr
The poset $\mathsf{s}(\cD(n))$ is a sublattice of $\bbN^r$.
\epr

\bprf
Consider hyperbinary expansions $c,d,d'\in\cD(n)$ such that $c\lt d$ and $c\lt d'$, and assume $d\ne d'$.
By Lemma~\ref{cov}, there exist indices $i,j$ such that $d$ is obtained from $c$ by replacing $c_i2$ with $(c_i+1)0$, and $d'$ is obtained from $c$ by replacing $c_j2$ with $(c_j+1)0$.
In order for both covers to be well defined, we must have $|i-j|\ge 2$.
Hence, we may perform both moves simultaneously to construct an element $e\in\cD(n)$ such that $d\lt e$ and $d'\lt e$.
Using Proposition~\ref{prop:coversum}
we have $\mathsf{s}(e)\gt\mathsf{s}(d),\mathsf{s}(d')$ in $\bbN^r$. Hence, $\mathsf{s}(e)$ is the join of $\mathsf{s}(d)$ and $\mathsf{s}(d')$ in $\bbN^r$.
By Lemma~\ref{lem:sublattice_local}, we conclude that $\mathsf{s}(\cD(n))$ is a join subsemilattice of $\bbN^r$.

Using a similar argument, we can show that $\mathsf{s}(\cD(n))$ is a meet subsemilattice of $\bbN^r$. 
Therefore, $\mathsf{s}(\cD(n))$ is a sublattice of $\bbN^r$.
\eprf

Since $\bbN^r$ is distributive, and distributivity is inherited by sublattices, the following corollary is immediate.

\bco
The poset $\cD(n)$ is a distributive lattice.
\eco

By the Fundamental Theorem of Finite Distributive Lattices, $\cD(n)$ is isomorphic to the lattice of order ideals of its subposet of join-irreducible elements.
The following lemma describes all join-irreducibles.

\ble\label{lem:joinirr}
Given $i\in[r]$, let $m,q$ be the unique nonnegative integers such that we have $n = q\cdot 2^{k-i} + m$ and $0\le m < 2^{k-i}$.
If $c=\zh(q),\ d=\zh(m)$ and $z=0^{k-\ell(c)-\ell(d)}$, then the word $czd$ is a join-irreducible element of $\cD(n)$.
Moreover, every join-irreducible is of this form.
\ele

\bprf
An element of a lattice is join-irreducible if and only if it covers a unique element. 
Thus, by Lemma~\ref{cov}, $b\in\cD(n)$ is join-irreducible if and only if there is a unique index $j$ such that $b_j\ne 0$ and $b_{j+1}=0$.

Fix an index $1\le i\le r$, and let $m$, $q$, $c$ and $d$ be as defined in the statement of the lemma.
We observe that $\ell(c)=i$ and $\ell(d)\le k-i$.
Hence, we may define the possibly empty word $z=0^{k-\ell(c)-\ell(d)}$ such that $czd$ has length $k$.
This word is a hyperbinary expansion of $n$ since
\begin{align*}
    s(czd) &= s(c)\cdot 2^{k-\ell(c)} + s(zd)\\
    &= q \cdot 2^{k-i} + m\\
    &= n.
\end{align*}

We must show that $czd$ has a unique index $j$ in $czd$ with $(czd)_j\neq 0$ and $(czd)_{j+1}=0$.
Now $q\neq 0$ since $n\ge 2^{k-1}$ and $2^{k-i}\le 2^{k-1}$.   So, by Corollary~\ref{cor:min}, $c=\zh(q)$  contains no such index but ends with a nonzero element.
On the other hand we can assume that $m<2^{k-i}-1$ since if $m=2^{k-i}-1$ then
$n=(q+1)2^{k-i}-1$ which implies that $r=0$. So in this case there is no $i\in[r]$.
Appealing to Corollary~\ref{cor:min} again, we see that $d=\zh(m)$ contains no index  as above and either starts with a zero element or is empty.  And in the latter case, $z$ will contain at least one zero.  Thus $j=\ell(c)$ is the unique index we seek in $czd$, making it join-irreducible.


Now suppose $e\in \cD(n)$ is join-irreducible, and let $i$ be the unique index such that $e_i\ne 0$ and $e_{i+1}=0$. 
We again define integers $m,q$ such that $n=q\cdot 2^{k-i}+m$ and $0\le m<2^{k-i}$.
As noted in the proof of Proposition~\ref{sEmb},  $e_j= 1$ for $j>r+1$.
Since $e_{i+1}=0$ we must have $i+1\le r+1$ which implies $i\in[r]$.
Consider subwords $c=e_1\ldots e_i$ and $d'=e_{i+1}\ldots e_k$. 
Since $n=s(c)\cdot 2^{k-i} + s(d')$, we deduce that 
$s(d')\equiv m\bmod{2^{k-i}}$.
Since $e_{i+1}=0$, we have
\[ s(d') \le \sum_{j=i+2}^k 2\cdot 2^{k-j} < 2\cdot 2^{k-i-1} = 2^{k-i}. \]
Hence, $s(d') = m$ and $s(c) = q$.
Since $e$ is join-irreducible, the $0$ entries in $c$ and $d'$ must occur at the beginning of each word.
Therefore, the words $c$ and $d'$ must correspond to minimum elements of $\cD(q)$ and $\cD(m)$, respectively, possibly with extra $0$ entries in between to ensure the length of $cd'$ is $k$.
\eprf

As an example, we construct the join-irreducible elements of $\cD(10)$ using Lemma~\ref{lem:joinirr}.
Taking $i=1$, we have $10=1\cdot 2^{4-1} + 2$, so $q=1,\ m=2$.
Since $\zh(1)=1$ and $\zh(2)=02$, the corresponding element of $\cD(10)$ is $1002$.
If $i=2$, then $q=2,\ m=2$, and the corresponding element of $\cD(10)$ is $0202$.
Finally, if $i=3$, then $q=5,\ m=0$, and the corresponding element of $\cD(10)$ is $0210$.
Observe that $1002,\ 0202,\ 0210$ are the three join-irreducible elements of $\cD(10)$ shown on the left side of Figure~\ref{smap}.

\bpr
For $e\in\cD(n)$, $e$ is a join irreducible element of $\cD(n)$ if and only if $\tilde{\mathsf{s}}(e)$ is the indicator vector of a principal order ideal of $\cF(n)$.
\epr

\bprf
Let $\beta(n)=b_1\ldots b_k$ be the binary expansion of $n$, and let $f=\zh(n)$.  We can assume that $n\neq 2^k-1$ since otherwise $\cD(n)$ has no join irreducibles.  So,
by Corollary~\ref{cor:min}, we have $f = 0(b_2+1)\ldots(b_r+1)21^{k-r-1}$.
Using the definition of $r$, we have
\[ f_l = 
\begin{cases}
    b_l-1\ &\mbox{if } l=1\\
    b_l+1\ &\mbox{if } 2\le l\le r\\
    b_l+2\ &\mbox{if } l=r+1\\
    b_l\ &\mbox{if } r+2\le l\le k
\end{cases}.
\]

Let $e$ be a join-irreducible element of $\cD(n)$. 
By Lemma~\ref{lem:joinirr}, there is a decomposition $e=czd$ and integers $i\in[r],q,m$ satisfying $n=q\cdot 2^{k-i}+m$, $0\le m<2^{k-i}$, such that $c=\zh(q),\ d=\zh(m)$, and $z=0^{k-\ell(c)-\ell(d)}$.

By definition, $q = s(b_1\ldots b_i)$ and $m=s(b_{i+1}\ldots b_k)$.
Since $i\in[r]$, the word $b_1\ldots b_i$ does not have any $2$'s.  And $b_1=1$ since it is the first digit of $\be(n)$.  So $b_1\ldots b_i$
is the binary expansion of $q$.
On the other hand, the word $b_{i+1}\ldots b_k$ may have some leading $0$'s, so the binary representation of $m$ may be a proper subword, say $\beta(m) = b_{j+1}\ldots b_k$ for some $j$ with $j\ge i+1$.  Note that if $m=0$ then $\be(m)$ is the empty word which is obtained by letting $j=k$.  Also, in the proof of Lemma~\ref{lem:joinirr} we showed that $m\neq 2^{k-i}-1$ so that, by definition of $r$ again, $j\le r+1$.
By the same token, if $j=r+1$ then $j+1=r+2$ and $b_{j+1}\ldots b_k=1^{k-j}$ while if $j<r+1$ then there is at least one zero in $b_{j+1}\ldots b_k$.
Since $e_{j+1} = d_1 = \zh(m)_1$, it follows from the previous sentence that either $j\le r$ and $e_{j+1}=0$, or $j=r+1$ and $e_{j+1}=1$. 
Similar reasoning shows that that $e_l=f_l$ for $l>j+1$.

We now prove that $\tilde{\mathsf{s}}(e)$ is the indicator vector of a principal order ideal of $\cF(n)$.
We first consider indices $l\le i$.
If $\beta(q)=1^i$, let $i'=1$.
Otherwise, let $i'\le i$ be the largest index such that $b_{i'}=0$.
In the former case, we have $c = 1^i$, whereas in the latter case, $c = 0(b_2+1)\ldots(b_{i'-1}+1)21^{i-i'}$.
Either way, for $1\le l\le i$,
\[ e_l - f_l = 
\begin{cases}
    0\ &\mbox{if }1\le l<i'\\
    1\ &\mbox{if }l=i'\\
    -1\ &\mbox{if }i'<l\le i
\end{cases}.
\]
Consequently, for $1\le l\le i$,
\[ s_l(e) - s_l(f) = 
\begin{cases}
    0\ &\mbox{if }1\le l<i'\\
    1\ &\mbox{if }i'\le l\le i
\end{cases}.
\]

Next consider $i<l\le\min\{j,r\}$.  We have $e_l=0$ and $f_l=1$ 
so that $e_l-f_l=-1$.  Now, reasoning as in the case $i'\le l\le i$ of the previous paragraph, we have
\beq
\label{ilmin}
 s_l(e) - s_l(f)  = 1. 
\eeq
Thus $\fst(e)$ is determined for this range of $l$.

For the remaining entries of $\tilde{\mathsf{s}}(e)$ when $l>\min\{j,r\}$ we separate two cases.

\vspace{3mm}

{\em Case 1: Assume $j\le r$.}

As stated above, we have $e_{j+1}=0$ in this case.  Also, $b_{j+1}=1$, since it is the leading digit of $\be(m)$, which implies $f_{j+1}=2$.
Also, in this case $\min\{j,r\}=j$ so that we have calculated $s_j(e)-s_j(f)=1$ in~\eqref{ilmin}. Hence
\[ s_{j+1}(e) - s_{j+1}(f) = (s_j(e) - s_j(f))\cdot 2 + (0-2) = 0. \]
If $l>j+1$, then $e_l=f_l$ by Corollary~\ref{cor:min}, which implies $s_l(e)=s_l(f)$ by induction.

\vspace{3mm}

{\em Case 2: Assume $j>r$. }
 
In this case, we have $j=r+1$
because of the restrictions placed on $m$ by the restriction $i\in[r]$ in Lemma~\ref{lem:joinirr}.
But $e_{r+1}=0$ and $f_{r+1}=2$, so we again find $s_{r+1}(e)-s_{r+1}(f) = 0$.
If $l>r+1$, then $e_l=f_l$, which implies $s_l(e)=s_l(f)$ by induction.

\vspace{3mm}

Let $I=\{x_l:i'\le l\le\min\{j,r\}\}$.
We have shown that $\tilde{\mathsf{s}}(e) = \chi_I$.  We claim that $I$ is the principal order ideal of $\cF(n)$ generated by $x_i$.
To prove this claim, we show that $I$ is the union of the principal order ideals generated by $x_i$ in each of the subposets $S = \{x_1,\ldots,x_i\}$ and $S'=\{x_i,\ldots,x_r\}$.

First, we show that $I\cap S=\{x_{i'},x_{i'+1},\ldots,x_i\}$ is the principal order ideal of $S$ generated by $x_i$. 
If $i'=1$, then $b_1=\ldots= b_i=1$, which implies $x_1\lt x_2\lt\cdots\lt x_i$.
If $i'\ge 2$, then by definition of $i'$, $b_{i'+1}=\ldots=b_i=1$ and $b_{i'}=0$.
Hence, $x_{i'}\lt x_{i'+1}\lt\cdots\lt x_i$ and $x_{i'-1}\gt x_{i'}$ in this case. 
In both cases, we conclude that $I\cap S$ is the principal order ideal of $S$ generated by $x_i$.

Next, we show that $I\cap S'$ is the principal order ideal of $S'$ generated by $x_i$.
We consider two cases. 

\vspace{3mm}

{\em Case 1: Assume $j\le r$.}
In this case, $I\cap S'=\{x_i,\ldots, x_j\}$.
By definition of $j$, $b_{i+1}=\cdots=b_j=0$ and $b_{j+1}=1$.
Then $x_i\gt x_{i+1}\gt\cdots\gt x_j$, and if $j<r$ then $x_j\lt x_{j+1}$.
Hence, the claim holds.

\vspace{3mm}

{\em Case 2: Assume $j>r$.}
In this case, $I\cap S'=\{x_i,\ldots,x_r\}$.
By definition of $j$, $b_{i+1}=\cdots=b_r=0$, which implies $x_i\gt x_{i+1}\gt\cdots\gt x_r$.
Since $x_r$ is the final element of the fence $\cF(n)$, we again conclude that the claim holds in this case. 

\vspace{3mm}

We have now completed the proof that $I$ is the principal order ideal generated by $x_i$.
Therefore, $\tilde{\mathsf{s}}(e)$ is the indicator vector of a principal order ideal of $\cF(n)$ whenever $e\in\cD(n)$ is join irreducible.
Since $|\cF(n)|=r$, and we have constructed $r$ join irreducible elements in Lemma~\ref{lem:joinirr}, the converse is true as well.
\eprf

Let $\irr(L)$ be the set of join-irreducible elements of a lattice $L$.

\ble
Let $S$ and $T$ be sublattices of a finite lattice $L$. If $\irr(S) = \irr(T)$, then $S=T$.
\ele

\bprf
We will just show $S\sbe T$ as the proof of the other containment is obtained by switching the roles of $S$ and $T$ in the demonstration. Let $x\in S$, and let 
$$
J = \{j\in\irr(S) : j\le x\}.
$$
Then $x=\bigvee_S J$.
Since $\irr(S)=\irr(T)$, we have the inclusion $J\subseteq T$.
And $S$ and $T$ are sublattices of the same lattice $L$ so that $\bigvee_S J = \bigvee_T J$.  Thus $x\in T$.
\eprf

We can now demonstrate the main result of this section.

\bth\label{thm:mainbij}
We have the poset isomorphism $\cD(n)\iso\cJ(\cF(n))$.
\eth
\bprf
We have shown that $\cD(n)$ and $\cJ(\cF(n))$ are each isomorphic to a sublattice of $\{0,1\}^r$, and these sublattices have the same set of join-irreducible elements.
Therefore, these two sublattices coincide.
\eprf

For $n=10$, the isomorphism of Theorem~\ref{thm:mainbij} is 
obtained by comparing the posets in Figures~\ref{fig:iso} and~\ref{smap}.

The {\em rank-generating function} of $\cJ(\cF(n))$ is 
\[ \rgf_n(t) = \sum_{I\in\cJ(\cF(n))} t^{|I|}. \]

\bth\label{thm:weightbij}
For all $n\in\bbN$, if the principal prefix of $\beta(n)$ has length $r$ and $\beta(n)$ has $s$ ones, then
\beq\label{eqn:htorgf}
 h_q(n) = q^{r+s} \rgf_n(q^{-1}).
\eeq
\eth

\bprf
We have 
\beq\label{eqn:rgf}
q^{r+s}\rgf_n(q^{-1}) = \sum_{I\in\cJ(\cF(n))} q^{r+s-|I|}.
\eeq
Applying the isomorphism in Theorem~\ref{thm:mainbij}, we claim that this sum transforms into $h_q(n)$.
By the assumptions of the theorem, the element $\be(n)=\oh(n)$  viewed as  a  partition in $\cD(n)$ has length $s$ and so contributes $q^s$ to $h_q(n)$.
On the other hand, $\oh(n)$ corresponds to the maximum ideal of $\cF(n)$ which contains all $r$ elements, so the analogous term 
in~\ref{eqn:rgf} is $q^{r+s-r}=q^s$.
By Lemma~\ref{cov}, removing an element from an ideal corresponds to increasing the length of the associated hyperbinary partition by $1$.
Therefore, Equation~\ref{eqn:htorgf} holds.
\eprf

Using the definition of $\CW_q(n)$, we deduce the following corollary from the previous theorem.

\bco\label{cor:qCW_fence}
Given $n\geq 1$, suppose $|\cF(n-1)|=r',\ |\cF(n)|=r,\ \beta(n-1)$ has $s'$ ones, and $\beta(n)$ has $s$ ones.
Then
$$
\hspace*{150pt} 
\CW_q(n) = q^{r'-r+s'-s} \, \frac{\rgf_{n-1}(q^{-1})}{\rgf_n(q^{-1})}. 
\hspace*{150pt}\qed
$$
\eco

\begin{figure}
\centering
\begin{tikzpicture}
\begin{scope}
    \draw (0,0) circle(1mm);
    \draw (1,0) circle(1mm);
    \draw (2,0) circle(1mm);
    \draw[<-] (.1,0) -- (.9,0);
    \draw[->] (1.1,0) -- (1.9,0);    
\end{scope}
\begin{scope}[xshift=-2cm, yshift=1cm]
    \filldraw (0,0) circle(1mm);
    \draw (1,0) circle(1mm);
    \draw (2,0) circle(1mm);
    \draw[<-] (.1,0) -- (.9,0);
    \draw[->] (1.1,0) -- (1.9,0);    
\end{scope}
\begin{scope}[xshift=2cm, yshift=1cm]
    \draw (0,0) circle(1mm);
    \draw (1,0) circle(1mm);
    \filldraw (2,0) circle(1mm);
    \draw[<-] (.1,0) -- (.9,0);
    \draw[->] (1.1,0) -- (1.9,0);    
\end{scope}
\begin{scope}[yshift=2cm]
    \filldraw (0,0) circle(1mm);
    \draw (1,0) circle(1mm);
    \filldraw (2,0) circle(1mm);
    \draw[<-] (.1,0) -- (.9,0);
    \draw[->] (1.1,0) -- (1.9,0);    
\end{scope}
\begin{scope}[yshift=3cm]
    \filldraw (0,0) circle(1mm);
    \filldraw (1,0) circle(1mm);
    \filldraw (2,0) circle(1mm);
    \draw[<-] (.1,0) -- (.9,0);
    \draw[->] (1.1,0) -- (1.9,0);    
\end{scope}
\begin{scope}[yshift=4cm]
    \draw(-1,0) node{$\cG = $};
    \filldraw (0,0) circle(1mm);
    \filldraw (1,0) circle(1mm);
    \filldraw (2,0) circle(1mm);
    \draw[<-] (.1,0) -- (.9,0);
    \draw[->] (1.1,0) -- (1.9,0);    
\end{scope}

\end{tikzpicture}
\caption{The directed graph $\cG$ for $r/s=5/2$ and its closure sets \label{fig:MGOgraph}}
\end{figure}
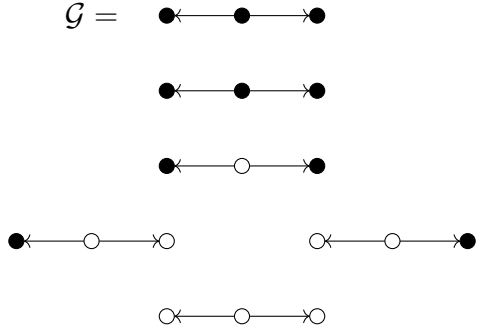

Corollary~\ref{cor:qCW_fence} may be alternately obtained from Theorem~\ref{thm:qrat} by using a different interpretation of $q$-deformed rationals discovered by Morier-Genoud and Ovsienko, as we now explain.
Suppose $r/s>1$ is a rational number with continued fraction expansion given by equation~\eqref{ConFra}.
Let $N=a_1+\cdots+a_m$ and consider the path graph with $N$ edges on a horizontal line.
In groups from left to right, we orient $a_1$ edges to the left, then $a_2$ edges to the right, then $a_3$ edges to the left, and so on.
Finally, we delete the two vertices on the ends to obtain a directed path graph $\cG$.
For example, the $11$th Calkin-Wilf number is $\CW(11)=5/2$, which has continued fraction representation $[2,2]$ and the corresponding $\cG$ is shown in the top line of Figure~\ref{fig:MGOgraph}.

A subset $X$ of vertices of $\cG$ is a {\em closure set} if there is no arrow $u\ra v$ such that $u\in X$ and $v\notin X$.
Let 
$$
f_{\cG}(q) = \sum_X q^{|X|}
$$ 
where the sum is over all closure sets $X$ for $\cG$. 
In Figure~\ref{fig:MGOgraph} the closure sets $X$ of $\cG$ are displayed in the last four lines, where a vertex is black or white depending on whether the vertex is or is not in $X$, respectively.
We also consider a subgraph $\cG'$ of $\cG$ obtained by deleting an additional $a_1$ vertices from the left side of $\cG$. With this setup, \cite[Theorem 4]{MGO:qdr1} 
(acknowledging ties with work of Lee and Schiffler~\cite{LS:cajp}) 
states that 
\beq
\label{GG'}
\left[\frac{r}{s}\right]_q = \frac{f_{\cG}(q)}{f_{\cG'}(q)}.
\eeq
Continuing our example, we see from  Figure~\ref{fig:MGOgraph} that $f_{\cG}(q) = 1+2q+q^2+q^3$.
The other graph $\cG'$ only has one vertex, and $f_{\cG'}(q) = 1+q$.
Hence, $[5/2]_q = \frac{1+2q+q^2+q^3}{1+q}$.

To compare~\eqref{GG'} with Corollary~\ref{cor:qCW_fence}, note that the directed graphs $\cG,\cG'$ can be converted into a fence poset where a directed edge $u\ra v$ is replaced by a cover $u\lt v$.
Under this correspondence, closure sets become order filters, which are in bijection with order ideals by complementation.  Continuing our example,  we have $\rgf_{10}(q) = 1+q+2q^2+q^3$ from Figure~\ref{cJcA10}.
Since $\cF(11)$ is a 1-element poset, its lattice of order ideals has rank generating function $\rgf_{11}(q)=1+q$.  Using Theorem~\ref{thm:qrat} and
then  Corollary~\ref{cor:qCW_fence}, we see that
\[
[5/2]_q = [\CW(11)]_q = q\CW_q(11) 
= q\cdot q\frac{\rgf_{10}(q^{-1})}{\rgf_{11}(q^{-1})} 
= \frac{1+2q+q^2+q^3}{1+q}
\]
which agrees with the previous computation.

\section{Matrices}
\label{m}

Consider the matrices 
\beq
\label{LR}
L = \left[
\barr{cc}
1 & 0\\
1 & q^{-1}
\earr
\right]
\qmq{and}
R = \left[
\barr{cc}
q & 1\\
0 & 1
\earr
\right].
\eeq
Mourier-Genoud and Ovsienko~\cite{MGO:qdr1} showed that the $q$-analogues of rational numbers can be expressed as ratios of entries in products involving $L$ and $R$.  
In this section, we will relate certain products to hyperbinary partitions.  We note that Han et al.~\cite{HEA:uvCW,HEA:mrv} have studied a generalization of the Calkin-Wilf sequence generated by matrices
$$
L_u = \left[
\barr{cc}
1 & 0\\
u & 1
\earr
\right]
\qmq{and}
R_v = \left[
\barr{cc}
1 & v\\
0 & 1
\earr
\right].
$$

Define a sequence of matrices $M(n)$ for $n\ge1$ as follows.
Let the binary expansion of $n$ be $\be(n)=b_1\ldots b_k$ so that $b_1=1$.  Removing the initial $1$ and reading the sequence backwards results in $b_k b_{k-1}\ldots b_2$.  Now let $M(n)$ be the matrix obtained from the product formed by replacing each $0$ in $b_k b_{k-1}\ldots b_2$ by $L$ and each $1$ by $R$.  For example, if $n=19$ then $\be(19)=10011$.  So, the reversed sequence is $1100$ and 
\begin{align*}
M(19)
&= RRLL\\
&= \left[
\barr{cc}
q & 1\\
0 & 1
\earr
\right] 
\left[
\barr{cc}
q & 1\\
0 & 1
\earr
\right]
\left[
\barr{cc}
1 & 0\\
1 & q^{-1}
\earr
\right]
\left[
\barr{cc}
1 & 0\\
1 & q^{-1}
\earr
\right]
\\[5pt]
&=
\left[
\barr{cc}
q^{-1}+ 2 + q + q^2 & q^{-1} + q^{-2}\\
q^{-1} + 1          & q^{-2}
\earr
\right].
\end{align*}


We give a formula for the entries of $M(n)$ in Theorem~\ref{thm:MnEnt}. 
First, we compute formulas for $h_q(n)$ for certain $n$ in the following lemma.

\ble\label{lem:bdcases}
Let $n$ be such that one has the factorization
$$
\be(n) = ac
$$
where $a=1^r$ is the principal prefix, and $c$ has at most one zero together with $s$ ones.  Then
$$
h_q(n)=\case{q^{r+s} + q^{r+s+1}+\cdots+q^{2r+s}}{if $c$ has a zero,}
{q^r}{if $c$ has no zero.}
$$
\ele
\bprf
Since $a$ is the principal prefix, the assumptions on $c$ imply that  either $c$ is the empty word or $c=01^s$.  In the former case, $a=1^r$ is the unique hyperbinary expansion of $n$
and has length $r$.  Thus $h_q(n)=q^r$.

If $c=01^s$ then, from Proposition~\ref{prop:bounds}, we have
$$
\oh(n) = \be(n)=1^r 0 1^s.
$$
Now Lemma~\ref{cov} implies that $\cD(n)$ is a chain with $r+1$ elements.  
Also, $\ell(\oh(n))=r+s$ so this element contributes $q^{r+s}$ to $h_q(n)$.
By the lemma just cited, length increases by $1$ as one goes from an element of $\cD(n)$ to an element it covers.  The formula given for $h_q(n)$ in this case now follows.
\eprf

For $k,l\in\bbN$, we use the notation
$$
[k,l) =\{k,k+1,\ldots,l-1\}.
$$
\bth\label{thm:MnEnt}
Suppose that $n\in[2^k, 2^{k+1}-1)$ and write $\be(n)=b_1 b_2\ldots,b_{k+1}$.  Let $j$ be the maximum index such that $b_1=\ldots=b_j=1$ and define $n'$ by
$$
\be(n') = 1 b_{j+2} b_{j+3} \ldots b_{k+1}.
$$
Then
$$
M(n) = 
\left[
\barr{cc}
q^{-k+2j-1}\ h_q(n' -1)    & q^{-k+1}\ h_q(n - 2^k - 1)\\[5pt]
q^{-k+2j-2}\ h_q(n')         & q^{-k}\ h_q(n - 2^k)      
\earr
\right].
$$
If $n=2^{k+1}-1$ then the same formula holds with the first column replaced by
$$
\left[
\barr{c}
q^k  \\[5pt]
0              
\earr
\right].
$$
\eth

\bprf

For each $n\ge 1$, let $N(n)$ be the matrix on the right-hand side of the above equation. We prove $M(n)=N(n)$ by induction on $k$ and divide the proof into four cases.

\vspace{3mm}

{\em Case 1: $n = 2^{k+1}-1$. }

In this case, $\beta(n) = 1^{k+1}$, so $M(n)=R^k$.
By induction, one can check that
\[ R^k = 
\begin{bmatrix}
    q^k & 1+q+\cdots+q^{k-1}\\
    0 & 1
\end{bmatrix}.
\]
Applying the identity $n-2^k = 2^k-1$, we have
\[ N(n) = 
\begin{bmatrix}
    q^k & q^{-k+1} h_q( 2^k-2 )\\
    0 & q^{-k} h_q(2^k - 1)
\end{bmatrix}.
\]
Now $\be(2k-1)=1^k$ and $\be(2^k-2)=1^{k-1}0$.  So,
by Lemma~\ref{lem:bdcases}, we deduce $M(n)=N(n)$.

\vspace{3mm}

{\em Case 2: $n = 2^{k+1} - 2$. }

The binary expansion of $n$ is $\beta(n)=1^k0$, so $M(n)=LR^{k-1}$.
From Case 1, we have
\[ M(n) = 
\begin{bmatrix}
    1 & 0\\
    1 & q^{-1}
\end{bmatrix}
\begin{bmatrix}
    q^{k-1} & 1 + q+\cdots + q^{k-2}\\
    0 & 1
\end{bmatrix} = 
\begin{bmatrix}
    q^{k-1} & 1 + q+\cdots + q^{k-2}\\
    q^{k-1} & q^{-1} + 1 + q+\cdots + q^{k-2}    
\end{bmatrix}.
\]
Since there are $k$ leading $1$'s in the binary expansion of $n$, we have $j=k$ and $n'=1$.
We also have $n-2^k = 2^k-2$.
Hence,
\[ N(n) = 
\begin{bmatrix}
    q^{k-1} h_q(0) & q^{-k+1} h_q(2^k-3)\\
    q^{k-2} h_q(1) & q^{-k} h_q(2^k-2)
\end{bmatrix}.
\]
We have $\be(2^k-3)=1^{k-2}01$ and $\be(2^k-2)$ from Case 1.
Applying Lemma~\ref{lem:bdcases}, we obtain $M(n)=N(n)$ again.

\vspace{3mm}

{\em Case 3: $n\in[2^k,2^{k+1}-2)$ and $n$ is odd.}

Let $m\in\bbN$ such that $n=2m+1$.
Then $\beta(n)=b_1\ldots b_k1$ and $\beta(m)=b_1\ldots b_k$.
Hence, $M(n) = R\cdot M(m)$.

Since $\beta(n)\ne 1^{k+1}$, the words $\beta(n)$ and $\beta(m)$ have the same number of leading $1$'s.
In particular, we have $n' = 2m' +1$.
Since $m\in[2^{k-1},2^k-1)$ we may apply the inductive hypothesis to get
\begin{align*}
M(n) &= 
\begin{bmatrix}
    q & 1\\
    0 & 1
\end{bmatrix}
\begin{bmatrix}
    q^{-k+2j}h_q(m'-1) & q^{-k+2} h_q(m-2^{k-1}-1)\\
    q^{-k+2j-1}h_q(m') & q^{-k+1} h_q(m-2^{k-1})
\end{bmatrix}\\[10pt]
&=
\begin{bmatrix}
    q^{-k+2j+1}h_q(m'-1) + q^{-k+2j-1}h_q(m') & q^{-k+3}h_q(m-2^{k-1}-1)+q^{-k+1}h_q(m-2^{k-1})\\
    q^{-k+2j-1}h_q(m') & q^{-k+1}h_q(m-2^{k-1})
\end{bmatrix}\\[10pt]
&=
\begin{bmatrix}
    q^{-k+2j-1}h_q(2m') & q^{-k+1} h_q(2m-2^k)\\
    q^{-k+2j-2}h_q(2m'+1) & q^{-k} h_q(2m-2^k+1)
\end{bmatrix}
\end{align*}
where the last equation follows by applying Proposition~\ref{hqrr}.
Using the substitutions $n=2m+1$ and $n'=2m'+1$, we conclude $M(n)=N(n)$.

\vspace{3mm}

{\em Case 4: $n\in[2^k,2^{k+1}-2)$ and $n$ is even.}

The proof in this case is very similar to the proof of Case 3, starting with $n=2m,\ n'=2m'$, and $M(n)=L\cdot M(m)$.
Due to the similarity to Case 3, we omit the proof.
\eprf


The row sums of $M(n)$ take a particularly nice form.
This result can be derived as a corollary of Theorem~\ref{thm:MnEnt}, but we give a simpler proof using the definition of $M(n)$ and the recurrence in Proposition~\ref{hqrr}.
\bth
\label{MnThm}
If $n\in[2^k,2^{k+1})$ then
$$
M(n) 
\left[
\barr{c}
1\\[5pt]
1
\earr
\right]
= 
\left[
\barr{c}
q^{-k}\ h_q(n-1)\\[5pt]
q^{-k-1}\ h_q(n)
\earr
\right].
$$
\eth
\bprf
As usual, we induct on $n$ where there are two cases depending on parity.  We will only do the even case as the odd case is similar.  Since $2^k$ is the largest power of $2$ less than or equal to $n$, we have $\be(n)=b_1\ldots b_{k+1}$ and  $\be(2n)= \be(n) 0$.  Transforming this into matrices we see that we have $M(2n) = L M(n)$.  It follows from induction and Proposition~\ref{hqrr} that
\begin{align*}
 M(2n)  \left[
\barr{c}
1\\[5pt]
1
\earr
\right]
&
=L M(n)
\left[
\barr{c}
1\\[5pt]
1
\earr
\right]
\\[5pt]
&=
\left[
\barr{cc}
1 & 0\\
1 & q^{-1}
\earr
\right]
\left[
\barr{c}
q^{-k}\ h_q(n-1)\\[5pt]
q^{-k-1}\ h_q(n)
\earr
\right]
\\[5pt]
&=
\left[
\barr{c}
q^{-k}\ h_q(n-1)\\[5pt]
q^{-k}\ h_q(n-1)+q^{-k-2}\ h_q(n)
\earr
\right]
\\[5pt]
&=
\left[
\barr{c}
q^{-k-1}\ h_q(2n-1)\\[5pt]
q^{-k-2}\ h_q(2n)
\earr
\right]
\end{align*}
as desired.
\eprf

We can generalize Theorem~\ref{MnThm} as follows.  Let $r,s$ be two indeterminates and define
$$
L' = \left[
\barr{cc}
1 & 0\\
r & s
\earr
\right]
\qmq{and}
R' = \left[
\barr{cc}
r & s\\
0 & 1
\earr
\right].
$$
Also let $M'(n)$ be the matrix obtained from the $L,R$ product for $M(n)$ by replacing each $L$ by $L'$ and each $R$ by $R'$.  As for the hyperbinary polynomials, given a hyperbinary expansion 
$d=d_1d_2\ldots d_k$ for $n$ we consider the statistics
$$
t(d) = \text{number of twos in $d$},
$$
and
$$
z(d) = \text{numbers of nonleading zeros in $d$},
$$
that is, the number of zeros to the right of the leftmost nonzero digit of $d$.  For example, if $n=34$ and $d=020010$ then $t(d)=1$ and $z(d)=3$.  Define the generating function
$$
H_{rs}(n) = \sum_{d\in\cD(n)} r^{t(d)} s^{z(d)}.
$$
The proof of the following result is much the same as the demonstrations of Proposition~\ref{hqrr}  and Theorem~\ref{MnThm} and so is omitted.
\bth
\label{M'nThm}
We have $h_{rs}(-1)=0$, $h_{rs}(0) = 1$, and for $n\ge1$
\begin{align}
h_{rs}(2n-1) &=  h_{rs}(n-1),\\
h_{rs}(2n) &= s h_{rs}(n) + r h_{rs}(n-1).
\end{align}
For all $n\ge0$

\vs{10pt}

\eqed{
M'(n) 
\left[
\barr{c}
1\\[5pt]
1
\earr
\right]
= 
\left[
\barr{c}
h_{rs}(n-1)\\[5pt]
h_{rs}(n)
\earr
\right].
}
\eth
Of course, we could combine all three statistic, $\ell(d)$, $t(d)$, and $z(d)$.  The reader should be able to supply the details.

\section{Comments and future work}

We gather here various ideas for future work and open problems.

\subsection{Negative numbers}

In the context of the $q$-deformed rationals, it's very natural to include negative rationals along with positive ones. One can devise two-sided Stern sequences 
$$\dots -1, 3, -2, 3, -1, 2, -1, 1, 0, 1, 1, 2, 1, 3, 2, 3, 1, \dots$$
and
$$\dots 1, -3, 2, -3, 1, -2, 1, -1, 0, 1, 1, 2, 1, 3, 2, 3, 1, \dots,$$
either of which will allow us to obtain every rational number (along with the honorary number $\infty$) as a quotient of successive terms, 
but is there a rationale for either of these artificial-seeming sequences?

Perhaps a clue comes from the 2-adic numbers. In this context, $-1$ can be represented through the left-infinite digit-sequence $\cdots 111$. Might such representations provide a hint?

\subsection{Lattices}

As illustrated at the beginning of Section~\ref{php}, the set of all partitions of $n$ do not form a lattice under refinement for $n\ge5$.  
But the subset of hyperbinary partitions does and, in fact, the lattice is distributive.  It would be interesting to identify other natural subposets of the full partition poset which are lattices and satisfy various lattice properties.

\subsection{Other statistics}

At the end of Section~\ref{m} we indicated how a couple of our results could be modified using two other statistics.  It would be interesting to see whether other theorems in this work have such analogues.  Also there are other statistics that could be studied.

As another example, given a hyperbinary partition $\eta$ we let
$$
p_i(\eta) = \text{number of parts of $\eta$ of multiplicity $i$}
$$
for $i = 1, 2$.  In Section~\ref{m} we used the notation $t(\eta)=p_2(\eta)$.
As an alternative description, suppose that $\eta$ is written out in terms of its digits in hyperbinary $d=d_1 d_2\ldots d_k$.  Then
$$
p_i(\eta) = \text{number of digits in $d$  equal to $i$}.
$$
Note that we have the relation
\beq
\ell(\eta) = p_1(\eta) + 2p_2(\eta).
\label{ellp1p2}
\eeq
The statistics $p_1$ and $p_2$, have been considered, respectively, by
Klav\v{z}ar, Milutinovi\'c, and Petr~\cite{KMP:sp}
and by Bates and Mansour~\cite{BM:qCW}.
As far as we know our statistic $\ell$ has not been studied before,
though is the one most closely related to the $q$-rationals of Morier-Genoud and Ovsienko.

Consider the generating function
$$
\hb_{s,t}(n) = \sum_{\eta\in H(n)} s^{p_1(\eta)} t^{p_2(\eta)}.
$$
Using~\eqref{ellp1p2} we see that setting $s=q$ and $t=q^2$ recovers our previously considered polynomial $h_q(n) = \sum_{\eta\in H(n)} q^{\ell(\eta)}$.  
The next result is derived from Proposition~\ref{Hrr} in much the same way as Proposition~\ref{hqrr}, so we omit the details. 
\bpr
We have $\hb_{s,t}(0) = 1$, and for $n\ge1$
\begin{align*}
\hb_{s,t}(2n-1) &=  \hb_{s,t}(n-1),\\
\hs{140pt} \hb_{s,t}(2n) &= \hb_{s,t}(n) + q^2 \hb_q(n-1).\hs{140pt} \qed
\end{align*}
\epr

A different statistic was studied by Dilcher, Ericksen, and Stolarsky~\cite{DE:hesp},\cite{DS:pass};
their statistic reduces each nonzero digit by 1 and interprets the result in binary.

\subsection{Chip firing}

Another perspective that might give rise to analogues of our theorems is the chip-firing perspective.
One can regard the hyperbinary expansion $d_1 d_2 \dots d_k$ of a number as configurations of chips on the natural numbers, with $d_i$ chips residing at the location $k-i$, 
where a chip-firing move replaces 2 chips at $m$ by 1 chip at $m+1$, where a site with 3 or more chips {\em must} fire
and a site with 2 chips {\em may} fire.
Although our write-up does not mention chips explicitly, 
we found the perspective a helpful source of intuition.
Other chip-games might have similar properties.
Richard Stanley (in private communication) proposes
a generalization in which no location can have more than $r$ chips, 
and each chips divides into $s$ chips when it moves one step to the right.

\bigskip

\noindent
{\sc Acknowledgments:} The authors acknowledge helpful suggestions from
Neil Calkin, Sophie Morier-Genoud, Valentin Ovsienko,
Bruce Reznick, Richard Stanley, and G\"unter Ziegler.

\nocite{*}
\bibliographystyle{alpha}

\end{document}